\IfFileExists{pspum-l.cls}{
  \documentclass{pspum-l}
}{
  \documentclass{amsproc}
  \typeout{WARNING: pspum-l.cls not found; using amsproc for drafting only.}
}

\usepackage[T1]{fontenc}
\usepackage{amsmath,amssymb,mathtools}
\usepackage{graphicx}
\usepackage{booktabs}

\usepackage[colorlinks=true,
            linkcolor=blue,
            citecolor=blue,
            urlcolor=blue]{hyperref}
\usepackage[nameinlink,capitalize,noabbrev]{cleveref}

\usepackage{lmodern}
\usepackage{amsmath, amsthm, amssymb, amsfonts}
\usepackage[normalem]{ulem}
\usepackage{hyperref}

\usepackage{mathrsfs}

\usepackage{verbatim} 
\usepackage{longtable}

\usepackage{mathtools}

\usepackage{tikz}
\usetikzlibrary{decorations.pathmorphing}
\tikzset{snake it/.style={decorate, decoration=snake}}

\usepackage{caption}

\usepackage{tikz-cd}
\usetikzlibrary{arrows}

\theoremstyle{plain}
\newtheorem{thm}{Theorem}[section]

\newtheorem{prop}[thm]{Proposition}
\newtheorem{conj}[thm]{Conjecture}

\theoremstyle{definition}
\newtheorem{defn}[thm]{Definition}

\theoremstyle{remark}
\newtheorem{rmk}[thm]{Remark}

\newcommand{\BC}{{\mathbb{C}}}

\newcommand{\BH}{{\mathbb{H}}}

\newcommand{\BP}{{\mathbb{P}}}
\newcommand{\BQ}{{\mathbb{Q}}}
\newcommand{\BR}{{\mathbb{R}}}

\newcommand{\BZ}{{\mathbb{Z}}}

\newcommand{\CC}{{\mathcal C}}

\newcommand{\CH}{{\mathcal H}}

\newcommand{\CK}{{\mathcal K}}
\newcommand{\CL}{{\mathcal L}}

\newcommand{\CO}{{\mathcal O}}
\newcommand{\CP}{{\mathcal P}}

\newcommand{\heart}{\ensuremath\heartsuit}

\DeclareFontFamily{OT1}{rsfs}{}
\DeclareFontShape{OT1}{rsfs}{n}{it}{<-> rsfs10}{}
\DeclareMathAlphabet{\curly}{OT1}{rsfs}{n}{it}

\usepackage{tikz}
\usepackage{lmodern}
\usetikzlibrary{decorations.pathmorphing}

\title[Cohomology of compactified Jacobians]{Cohomology of compactified Jacobians for locally planar integral curves}

\author[J. Shen]{Junliang Shen}
\address{Yale University}
\email{junliang.shen@yale.edu}

\begin{document}

\begin{abstract}
This article surveys some recent developments on the cohomology of the compactified Jacobian associated with a locally planar integral curve. Topics discussed here include the Ng\^o support theorem, the perverse filtration, connections to the Hilbert schemes, and cohomological structures induced by the Arinkin--Fourier--Mukai transform.
\end{abstract}

\maketitle

% Optional table of contents for a long expository article.
% \tableofcontents

%\section{Introduction}

\setcounter{section}{-1}

\section{Introduction}

Throughout, we work over the complex numbers $\BC$. 

For a nonsingular irreducible projective curve $C$ of genus $g$, the space of holomorphic differential forms is a $g$-dimensional $\BC$-vector space $H^0(C, \Omega^1_C)$. The homology lattice $H_1(C,\BZ)$ embeds naturally in the dual vector space $H^0(C, \Omega^1_C)^\ast$ via integration
\[
[\gamma] \in H_1(C,\BZ) \mapsto \left(\omega  \mapsto \int_{\gamma} \omega \right) \in H^0(C, \Omega_C^1)^\ast.
\]
The associated quotient $H^0(C, \Omega^1_C)^*/H_1(C, \BZ)$ is a $g$-dimensional complex torus, called the Jacobian $J_C$ of the curve $C$. 

The Jacobian is naturally an abelian variety --- a complex torus which can be realized as an algebraic variety. An algebraic interpretation of $J_C$ is given by the Picard variety parameterizing degree 0 line bundles on the curve $C$,
\[
J_C = \mathrm{Pic}^0(C):=\{\mathrm{degree~0~ line~ bundles~ on~ C}\}.
\]
This modular description is often regarded as the definition of the Jacobian $J_C$ in algebraic geometry.

For a singular curve $C$, the Jacobian variety $J_C$ (parameterizing degree 0 line bundles) may no longer be proper.\footnote{The variety $J_C$ for a singular curve $C$ is sometimes called the generalized Jacobian.} A typical example is given by a nodal plane cubic $C$, where $J_C \simeq \BC^*$. If the curve $C$ is integral (\emph{i.e.} reduced and irreducible), there is a natural modular compactification of $J_C$ --- by \cite{AK} there is a fine projective moduli space $\overline{J}_C$ parameterizing torsion free sheaves of rank 1 and degree 0 on the integral curve $C$; this is called the compactified Jacobian associated with $C$. Furthermore, we know from \cite{AIK, KK} that $\overline{J}_C$ is irreducible if and only if the integral curve $C$ is locally planar (\emph{i.e.} each analytic neighborhood of a singular point of $C$ can be embedded in $\mathbb{C}^2$). In this case, the compactified Jacobian 
\[
J_C \subset \overline{J}_C
\]
is an ``honest'' compactification of $J_C$: it contains $J_C$ as a Zariski open and dense subset, and the boundary is given by those rank 1 torsion free sheaves that fail to be line bundles.

In this survey article, we focus on the compactified Jacobian associated with a locally planar integral curve $C$. There are (at least) three reasons for such a choice. First, it is a fundamental and simple geometric object, \emph{e.g.} it is integral and  \emph{l.c.i.}\footnote{Recall that a morphism is \emph{l.c.i.} if it can be factored locally as a regular embedding followed by a smooth morphism; a scheme is \emph{l.c.i.} if the natural morphism to a point is \emph{l.c.i}.}, and there is no complicated stability condition involved as in other moduli problems in algebraic geometry. Second, it carries rich structures that lie in the crossroads of many branches of mathematics. Third, many new tools have been developed in recent years to study them, \emph{e.g.} perverse sheaves, derived categories, representation theory, knot/link theory, \emph{etc}.

Let me give some examples of the second point above to illustrate that the cohomology of the compactified Jacobian associated with a locally planar integral curve has rich connections:

\begin{enumerate}
    \item[$\bullet$] Many compactified Jacobians are geometric models of affine Springer fibers, whose cohomology admits interesting representations \cite{L, OY, OY2, GORS}. 
    \item[$\bullet$] The topology of a planar singular point is characterized by the associated link; the cohomology of the compactified Jacobian of a locally planar curve is expected to detect deep link invariants of the singularities \cite{ORS, Sh0}; 
    \item[$\bullet$] Compactified Jacobians associated with spectral curves (which are automatically locally planar) arise as Hitchin fibers; therefore, the local cohomological study of the Hitchin fibration is reduced to the study of the cohomology of certain compactified Jacobians \cite{Ngo, CL, MS_chi};
    \item[$\bullet$]  Cohomology of compactified Jacobians is closely related to the Gopakumar--Vafa invariants in enumerative geometry \cite{MT,MY, MS, DK};
    \item[$\bullet$]  Compactified Jacobians are natural geometric degenerations of abelian varieties; therefore they can be used to test conjectures for abelian fibrations with singular fibers \cite{MSY, MSY2, BMSY1}.
\end{enumerate}

An important structure that is related to all the stories above is the \emph{perverse filtration}. In short, the cohomology $H^*(\overline{J}_C, \BQ)$ is not only a graded vector space, but is naturally \emph{filtered}:
\begin{equation}\label{perv_filtration}
P_0H^*(\overline{J}_C, \BQ) \subset P_1H^*(\overline{J}_C, \BQ) \subset  P_2H^*(\overline{J}_C, \BQ) \subset \cdots \subset H^*(\overline{J}_C, \BQ).
\end{equation}
This filtration endows the cohomology with a ``hidden grading'', which plays a crucial role in the connections to geometric representation theory, link invariants, Gopakumar--Vafa invariants, \emph{etc}.

The perverse filtration for a general projective morphism was used in the seminal paper of de Cataldo--Migliorini \cite{dCM0} in their proof of the decomposition theorem \cite{BBDG} via classical Hodge theory. Then de Cataldo--Hausel--Migliorini \cite{dCHM} discovered that the perverse filtration associated with the Hitchin fibration plays a crucial role in non-abelian Hodge theory, and formulated the (now proven) $P=W$ conjecture \cite{MS_PW, HMMS, MSY, Hoskins}. The perverse filtration (\ref{perv_filtration}) for the compactified Jacobian is viewed as a ``local analog'' of the perverse filtration for the Hitchin fibration, which seems to be more mysterious. For example, the cohomology of the moduli space of stable Higgs bundles is generated by \emph{tautological classes} \cite{Markman}; as a key step in all the known proofs of the $P=W$ conjecture, the perverse filtration of the Hitchin fibration can be completely described in terms of the tautological classes. However, for the compactified Jacobian associated with a general locally planar integral curve, it is not clear how to present all classes in the cohomology explicitly; therefore it is more difficult to describe the perverse filtration. Furthermore, the non-abelian Hodge theory for the moduli space of stable Higgs bundles is well-understood; motivated by connections to link theory, Shende \cite{Sh0} proposed that the perverse filtration (\ref{perv_filtration}) should also be matched with the weight filtration of some \emph{Betti moduli space} under a version of non-abelian Hodge theory; see also \cite[Conjecture 6.16 and Remark 6.17]{GKS}. Although this proposal has been studied intensively after its appearance (see \emph{e.g.} \cite{BBAMY, Trinh2, Trinh}), it is still mysterious what is the correct geometric formulation of the non-abelian Hodge theory for an arbitrary planar singularity.

These notes present several topics on recent developments in the cohomological study of compactified Jacobians. We focus on the particularly accessible yet rich case of locally planar integral curves. The goal is to illustrate some new ideas and techniques emerging in this area, in which the perverse filtration plays a central role.

\subsection*{Plan of the following sections}
 Section \ref{sec1} focuses on the support theorem of Ng\^o; the support theorem is a fundamental new input in the modern cohomological study of compactified Jacobians. As a consequence, the cohomology of the compactified Jacobian carries a \emph{canonical} perverse filtration, introduced by Maulik--Yun; we discuss this filtration, its properties, and some open questions in Section \ref{sec2}. In Section \ref{sec3}, we discuss the cohomological relation between compactified Jacobians and Hilbert schemes; this can be regarded as a ``toy model'' for the correspondence of Gopakumar--Vafa invariants and Pandharipande--Thomas invariants, where the perverse filtration plays a crucial role. The last section discusses the Arinkin--Fourier--Mukai transform for the compactified Jacobian which extends the classical Fourier--Mukai transform for abelian varieties; we explain how this duality for the derived category of coherent sheaves yields new structural results for the perverse filtration associated with the compactified Jacobian.

\subsection*{Acknowledgements}

These notes are closely related to the topics discussed by the author at the graduate student and postdoc bootcamp associated with the 2025 Algebraic Geometry Summer Research Institute at Colorado State University. The author is grateful to the organizers of the bootcamp for their efforts in making it possible. Special thanks are extended to all participants in the author’s mentored group: Younghan Bae, Alessio Cela, Bochao Kong, Andrés Ibáñez Núñez, Anibal Aravena, Francesca Rizzo, Jacob Cleveland, Shengxuan Liu, Yifan Wu, and Yifan Zhao. These notes are written for the Bootcamp Proceedings. Finally, the author would like to thank Davesh Maulik, Qizheng Yin, Yifan Zhao, and the referee for reading these notes carefully.

J.S.~was supported by the NSF grant DMS-2301474 and a Sloan Research Fellowship.

\section{Ng\^o support theorem}\label{sec1}

\subsection{Compactified Jacobian fibrations}\label{sec1.1}
One of the first constructions of compactified Jacobians by Igusa \cite{Igusa} viewed the compactified Jacobian for a singular curve as the limit of a family of Jacobian varieties for nonsingular curves. The purpose of this section is to explain that the cohomology of the compactified Jacobian is also the ``limit'' of the cohomology of Jacobians for nonsingular curves. A precise formulation requires the decomposition theorem of Beilinson, Bernstein, Deligne, and Gabber \cite{BBDG} and its refinement due to Ng\^o \cite{Ngo}. 

To state the main result (Theorem \ref{thm1.1}) of this section, we first introduce some notation. We fix $C$ to be a locally planar integral curve of arithmetic genus $g$, and we include it in a ``large enough'' family as follows. We let $p: \CC \to B$ be a flat family of integral projective curves of arithmetic genus $g$ with planar singularities over a nonsingular base $B$, where $C \subset \CC$ is the fiber over a point $0\in B$. We denote by $\pi: \overline{J}_\CC \to B$ the compactified Jacobian fibration, and assume that the total space~$\overline{J}_\CC$ is nonsingular. In particular, $\pi: \overline{J}_\CC \to B$ is flat with integral fibers of dimension $g$. The family $\pi: \overline{J}_{\CC} \to B$ contains $\overline{J}_C = \overline{J}_{\CC_0}$ as a fiber over $0\in B$, and a general fiber is the Jacobian variety $J_{\CC_b}$ associated with a nonsingular curve $\CC_b$. By \cite[Proposition 3.5]{MY}, such a family always exists.

We denote by $p^\circ: \CC^\circ \to B^\circ$ the locus of nonsingular curves with $\pi^\circ: J_{\CC^\circ} \to B^\circ$ the smooth Jacobian fibration. For any nonsingular curve $\CC_b$ with $b \in B^\circ$, its total cohomology is the exterior algebra of the first cohomology of the curve,
\[
H^*(J_{\CC_b}, \BQ) = \wedge^* H^1(\CC_b, \BQ).
\]
A relative version of this result is given by Deligne's decomposition theorem for smooth projective maps:
\begin{equation}\label{DT_smooth}
R\pi^\circ_* \BQ_{J_{\CC^\circ}} \simeq \bigoplus_{i=0}^{2g} \wedge^i R^1p^\circ_* \BQ_{{\CC^\circ}} [-i] \in D^b_c(B^\circ).
\end{equation}
Here every term in the right-hand side is a (shifted) local system, and the isomorphism holds in the bounded derived category $D^b_c(B^\circ)$ of constructible sheaves on $B^\circ$.

The following theorem extends (\ref{DT_smooth}) to all the singular fibers via the intermediate extension of perverse sheaves. It is a direct application of Ng\^{o}'s support theorem \cite{Ngo}; see \cite[Theorem 2.4]{MY}.

\begin{thm}\label{thm1.1}
We have
\begin{equation}\label{DT_main}
R\pi_* \BQ_{\overline{J}_\CC} \simeq \bigoplus_{i=0}^{2g} \mathrm{IC}_B\left(  \wedge^i R^1p_* \BQ_{{\CC^\circ}}\right)[-i -\dim B] \in D^b_c(B).
    \end{equation}
\end{thm}

Here $\mathrm{IC}_B(\CL)$ stands for the intermediate extension of a local system $\CL$ defined over a Zariski open subset of $B$. Since conventionally the intermediate extension of a local system is a perverse sheaf, the extra shift of $-\dim B$ is to make the indices still lie in $[0, 2g]$.

The two sides of (\ref{DT_main}) are of different flavors: the left-hand side collects cohomology of all the fibers $\overline{J}_{\CC_b}$; the right-hand side only depends on the nonsingular curves $p^\circ: \CC^\circ \to B^\circ$. In other words, it shows that the cohomology $H^*(\overline{J}_C, \BQ)$ is completely determined by the collection of $H^1(\CC_b, \BQ)$ for the nonsingular curves.

\subsection{Decomposition Theorem}\label{sec1.2}

Let $f: X \to B$ be a proper morphism between nonsingular varieties with equidimensional fibers. By the decomposition theorem \cite{BBDG}, we have
\begin{equation}\label{DT_general}
Rf_* \BQ_X \simeq \bigoplus_{i} {^\mathfrak{p}}\CH^i(Rf_* \BQ_X)[-i] \in D^b_c(B).
\end{equation}
Here ${^\mathfrak{p}}\CH^i(-)$ stands for the perverse cohomology --- taking the cohomology with respect to the perverse $t$-structure:
\[
{^\mathfrak{p}}\CH^i(-): D^b_c(B) \to \mathrm{Perv}(B).
\]
A deep result of \cite{BBDG} is that each perverse cohomology ${^\mathfrak{p}}\CH^i(Rf_* \BQ_X)$ on the right-hand side is a semisimple perverse sheaf. Recall that a simple perverse sheaf is of the form $\mathrm{IC}_Z(\CL)$ where 
\begin{enumerate}
    \item[(a)]  $Z$ is an irreducible closed subset of $B$,
    \item[(b)] $\CL$ is an irreducible local system defined on an open subset $U$ in the smooth locus of $Z$, and
    \item[(c)] $\mathrm{IC}_Z(-)$ stands for the composition of the intermediate extension from $U$ to $Z$, and the extension by $0$ from $Z$ to $B$.
\end{enumerate}
The closed subset $Z$ is called the \emph{support} of the simple perverse sheaf $\mathrm{IC}_Z(\CL)$. By definition, a semisimple perverse sheaf is a direct sum of simple perverse sheaves. An important invariant of the decomposition theorem for $f: X \to B$ is the collection of the {supports} given by all the simple summands that arise in the right-hand side of (\ref{DT_general}). 

Now we revisit Theorem \ref{thm1.1} from the perspective of supports. As a general result following from the decomposition theorem, the right-hand side of (\ref{DT_main}) must appear in the decomposition of the left-hand side. In fact, $R\pi_* \BQ_{\overline{J}_\CC}$ is semisimple, and the intermediate extensions of the local systems $\wedge^i R^1p^\circ_* \BQ_{{\CC^\circ}}$ have to appear by (\ref{DT_smooth}). In other words, the terms in the right-hand side are exactly all the summands with full support $B$, and Theorem \ref{thm1.1} is equivalent to saying that there is no summand with smaller support.

In general, it is very difficult to calculate all the supports for a proper morphism. For a large class of abelian fibrations (a proper map whose general fiber is an abelian variety), Ng\^o developed a framework in his proof of the fundamental lemma of the Langlands program \cite{Ngo} which obtains a strong control of the supports. In the following, we discuss the ideas of Ng\^o and explain the proof of Theorem \ref{thm1.1}.

\subsection{Ng\^{o}'s support theorem}\label{sec1.3}

The idea of Ng\^o was built on the Goresky--MacPherson inequality (see \cite[Section 7]{Ngo}) which we briefly review. Let $f: X\to B$ be a proper morphism as in Section \ref{sec1.2} where $d$ is the relative dimension.

\begin{prop}[Goresky--MacPherson inequality]\label{prop1.2}
Let $Z \subset B$ be a support of the decomposition of $Rf_*\BQ_X$. Then the following hold.
\begin{enumerate}
    \item[(a)] We have
    \[
\mathrm{codim}_B(Z) \leq d.
\]
\item[(b)] If a support $Z$ reaches equality in (a), then $Z$ can be detected by the constructible sheaf $R^{2d}f_* \BQ_X$; more precisely, there is an open $U\subset B$ with $U\cap Z\neq \emptyset$ such that a direct sum component of the sheaf $R^{2d}f_*\BQ_X|_U$ is supported on $U\cap Z$.
\end{enumerate}

\end{prop}

We note that the sheaf $R^{2d}f_*\BQ_X$ which appears in Proposition \ref{prop1.2} (b) calculates the fundamental classes of the irreducible components of the fibers.

The spirit of the proposition is that a support cannot be ``too small''. The proof relies on a clever use of the Verdier dualizing functor, combined with the fact that $Rf_*\BQ_X[\dim X]$ is self-dual, and is concentrated in the degrees in $[-\dim X, -\dim X+2d]$.

For example, the proposition implies that the decomposition theorem associated with a family of curves can have support at worst at codimension 1 in the base. However, the larger the fiber dimension $d$ is, the less powerful the proposition is.

Ng\^o observed that if the fibration $f: X \to B$ admits a fiberwise symmetry given by a commutative group scheme, then the inequality of Proposition \ref{prop1.2} can be greatly improved. Ng\^o's framework is for $\delta$-regular \emph{weak abelian fibrations}; for convenience, we focus on the compactified Jacobian fibration $\pi: \overline{J}_\CC \to B$ of Section \ref{sec1.1}.

Recall that the relative Jacobian 
\[
J_\CC :=\{(C_b, L) ~\Big{|} ~ L \in \mathrm{Pic}^0(\CC_b), ~ b \in B\}
\]
is a group scheme over $B$; it is a Zariski open subset of the relative compactified Jacobian
\[
J_\CC \subset \overline{J}_\CC.
\]
Moreover, this group scheme acts on $\overline{J}_\CC$ via tensor product:
\begin{equation}\label{fiber_sym}
a: J_\CC \times_B \overline{J}_\CC \to \overline{J}_\CC, \quad \left((C_b,L), (C_b,F)\right) \mapsto (C_b, F\otimes L).
\end{equation}

Now for any point $b\in B$, we let $\widetilde{\CC}_b \to \CC_b$ be the normalization, where the genus of the nonsingular curve $\widetilde{\CC}_b$ is the geometric genus of the possibly singular $\CC_b$. The commutative group $J_{\CC_b}$ can be expressed as an extension
\[
1 \to R_b \to J_{\CC_b} \to J_{\widetilde{\CC_b}} \to 1,
\]
known as the Chevalley decomposition. The map $J_{\CC_b} \to J_{\widetilde{\CC}_b}$ is given by the pullback of degree 0 line bundles along the normalization map $\widetilde{\CC}_b \to \CC_b$, and the first factor $R_b$ is an affine group scheme. The dimension of $R_b$, which is simply the difference of the arithmetic and the geometric genera of $\CC_b$, defines a function on the base
\[
\delta: B \to \BZ_{\geq 0}, \quad  b \mapsto \dim R_b.
\]

Let $Z \subset B$ be an irreducible closed subset. We define the $\delta$-value of $Z$ as $\delta(b)$ with $b$ a very general point in $Z$; equivalently, by semicontinuity we have
\[
\delta_Z:=\mathrm{min}\{\dim R_b ~\Big{|}~ b\in Z\}.
\]
Using the fiberwise symmetry (\ref{fiber_sym}) and a purity argument, Ng\^o observed that the Goresky--MacPherson inequality can be further ``squeezed''. In the case of the compactified Jacobian fibration $\pi: \overline{J}_{\CC} \to B$, Ng\^o's support theorem provides an enhancement of Proposition \ref{prop1.2}.

\begin{thm}[Ng\^o support theorem \cite{Ngo}]\label{thm1.3}
    Let $\pi: \overline{J}_\CC \to B$ be the compactified Jacobian fibration as above, and let $Z\subset B$ be a support of the decomposition of $R\pi_*\BQ_{\overline{J}_\CC}$. 
    \begin{enumerate}
        \item[(a)]  We have
        \[
        \mathrm{codim}_B(Z) \leq \delta_Z.
        \]
        \item[(b)] If a support $Z$ reaches equality in (a), then $Z$ can be detected by the constructible sheaf $R^{2g}\pi_*\BQ_{\overline{J}_\CC}$.
    \end{enumerate}
\end{thm}

Here the part (b) is in the same sense as Proposition \ref{prop1.2} (b).

We conclude this section by explaining that Theorem \ref{thm1.1} is an immediate consequence of Theorem \ref{thm1.3}.

\begin{proof}[Proof sketch of Theorem \ref{thm1.1}]
 Since the total space $\overline{J}_\CC$ is nonsingular, we have the classical Severi inequality, which says any irreducible subset $Z\subset B$ satisfies
 \begin{equation}\label{severi}
 \mathrm{codim}_B(Z) \geq \delta_Z.
 \end{equation}
Intuitively, the smoothness of $\overline{J}_\CC$ guarantees that the $\delta$-value of $Z$ cannot be ``too large'' compared to its codimension --- the reason is that the base $B$, viewed as a slice in the moduli stack of genus $g$ curves, has to meet transversally to every $\delta$-stratum; see \cite[Theorem 7.3]{CL} or \cite[Lemma 4.1]{MS_chi}.

Now we assume that $Z \subset B$ is a support. Since (\ref{severi}) holds for any $Z$, we obtain from Theorem \ref{thm1.3} (a) that the support $Z$ has to be of codimension $\delta_Z$. Then by the part (b), it can be detected by the sheaf $R^{2g}\pi_* \BQ_{\overline{J}_\CC}$. Since every fiber $\overline{J}_{\CC_b}$ is irreducible, we know that $H^{2g}(\overline{J}_{\CC_b}, \BQ)$ is 1-dimensional; in particular $R^{2g}\pi_* \BQ_{\overline{J}_\CC}$ is a rank 1 local system on the base $B$.\footnote{In fact, this is a trivial local system; but it is not needed in the proof of the theorem.} This concludes that $Z = B$.
\end{proof}

\subsection{Further remarks}

\begin{enumerate}
    \item[(1)] Migliorini--Shende \cite{MS} provided different proofs of Theorem \ref{thm1.1} using the geometry of curves. More precisely, if $Z$ is a support, Migliorini--Shende showed that there must be a ``nice'' curve lying in $Z$ by deformation theory. Then one can calculate essentially both sides of (\ref{DT_main}) to ensure that there is no room for an extra factor. Their methods also proved a more general statement concerning the Hilbert schemes. We will discuss the connection to Hilbert schemes more systematically in Section \ref{sec3}.

    \item[(2)] Various versions of support theorem have been proven for compactified Jacobians associated with non-integral curves. Migliorini--Shende--Viviani \cite{MSV} generalized Theorem \ref{thm1.1} to locally planar reduced curves; for reducible curves stability conditions are needed to construct good models of compactified Jacobians. Chaudouard--Laumon \cite{CL} proved support theorem for certain compactified Jacobian fibrations associated with non-reduced curves that arise from meromorphic Higgs bundles. Maulik--Shen \cite{MS_chi} further generalized the Chaudouard--Laumon argument and proved support theorem without the nonsingular assumption of the total space. See also \cite{dCHeM, MM, MMP} for a (definitely incomplete) list of references in this direction.
\end{enumerate}

\section{Perverse filtrations}\label{sec2}

\subsection{Definition}\label{sec2.1}

Let $C$ be a locally planar integral curve of arithmetic genus $g$. The perverse filtration on $H^*(\overline{J}_{C}, \BQ)$ can be thought of as the ``limit'' of the decomposition
\begin{equation}\label{decomp1}
H^*({J}_{\CC_b}, \BQ) = \bigoplus_{i=0}^{2g} \wedge^i H^1(\CC_b, \BQ)
\end{equation}
for nonsingular curves $\CC_b$, which we explain in the following.

As before, we put $C = \CC_0$ in a family $\CC \to B$ as in Section \ref{sec1}. Then the decomposition (\ref{decomp1}) for nonsingular curves $\CC_b \subset \CC^\circ$ specializes \emph{na\"ively} to the cohomological decomposition
\begin{equation}\label{boring}
H^*(\overline{J}_C, \BQ) = \bigoplus_{i=0}^{2g} H^i(\overline{J}_C, \BQ).
\end{equation}
Obviously, this is a boring specialization. Alternatively, the local system formed by 
\[
H^i(J_{\CC_b}, \BQ) = \wedge^i H^1(\CC_b, \BQ),\quad b\in B^\circ
\]
can be extended over the singular fibers via intermediate extension as in the right-hand side of Theorem \ref{thm1.1}; if we take the stalk for both sides over $0\in B$, we obtain in an interesting way of decomposing $H^*(\overline{J}_C, \BQ)$ into $2g+1$ pieces:
\begin{equation}\label{decomp2}
H^*(\overline{J}_C, \BQ) \simeq \bigoplus_{i=0}^{2g} Q_i,
\end{equation}
where each vector space $Q_i$ is the stalk of the summand $\mathrm{IC}_B\left(  \wedge^i R^1p_* \BQ_{{\CC^\circ}}\right)[-i -\dim B]$ on the right-hand side of (\ref{DT_main}) viewed as a graded $\BQ$-vector space. This decomposition turns out to be highly different from the cohomological decomposition (\ref{boring}), which refines the cohomology group with a second mysterious grading other than the cohomological one.

\begin{defn}\label{def2.1}
The perverse filtration   
  \begin{equation}\label{P_fil}
P_0H^*(\overline{J}_C, \BQ) \subset P_1H^*(\overline{J}_C, \BQ) \subset  P_2H^*(\overline{J}_C, \BQ) \subset \cdots \subset P_{2g}H^*(\overline{J}_C, \BQ) \subset H^*(\overline{J}_C, \BQ)
\end{equation}
is an increasing filtration given by 
\[
P_k H^*(\overline{J}_C, \BQ):= \bigoplus_{i\leq k} Q_i \subset H^*(\overline{J}_C, \BQ).
\]
\end{defn}

Although (\ref{DT_main}) gives a decomposition, we note that it is not canonical, therefore it does not induce a canonical decomposition of $H^*(\overline{J}_C, \BQ)$. Instead, the filtration defined above is always canonical, as the sheaf-theoretic filtration
\[
\bigoplus_{i=0}^{k}  \mathrm{IC}_B\left(  \wedge^i R^1p^\circ_* \BQ_{{\CC^\circ}}\right)[-i -\dim B] \to R\pi_*\BQ_{\overline{J}_\CC}
\]
is induced by the perverse truncation functor
\[
{^\mathfrak{p}}\tau_{\leq k+\dim B} R\pi_* \BQ_{\overline{J}_\CC} \to R\pi_*\BQ_{\overline{J}_\CC}.
\]

\begin{rmk}
We will discuss other possible interesting splittings in Section \ref{decomp_?}. 
\end{rmk}

Now we have defined a filtration on $H^*(\overline{J}_C, \BQ)$ dependent on a choice of embedding $C$ into a family $\CC \to B$. In the next section, we show that the perverse filtration in fact is independent of that choice --- in particular, it defines an intrinsic structure on the cohomology $H^*(\overline{J}_C, \BQ)$.

We conclude this section by explaining why the decompositions (\ref{decomp1}) and (\ref{decomp2}) are different in general. Since $\overline{J}_C$ is singular, the cohomology $H^*(\overline{J}_C, \BQ)$ may fail the numerical Poincar\'e duality, \emph{i.e.}, there exists $i\in \BZ_{\geq 0}$ such that 
\[
\dim H^i(\overline{J}_C, \BQ) \neq \dim H^{2g-i}(\overline{J}_C, \BQ).
\]
On the other hand, by the relative Hard Lefschetz, we always have 
\[
\dim \mathrm{Gr}^P_i H^*(\overline{J}_C, \BQ) = \dim \mathrm{Gr}^P_{2g-i} H^*(\overline{J}_C, \BQ).
\]

\subsection{Examples}\label{Examples}
We discuss the geometry of compactified Jacobians and perverse filtrations through some examples.

An interesting class of compactified Jacobians are given by the ones associated with a rational curve with a singular point given by the equation
\begin{equation}\label{p_q}
x^p-y^q=0, \quad \mathrm{gcd}(p,q)=1.
\end{equation}
These compactified Jacobians have deep connections to geometric representation theory and knot theory; see \cite{GORS, GKS, KR, OY,OY2}.

We first consider the case $(p,q)=(2,3)$. It is clear that $C$ is a cubic cuspidal curve which is isomorphic to $\overline{J}_C$. Therefore the compactified Jacobian admits an affine paving
\[
\overline{J}_C = \BC^1 \sqcup \BC^0.
\]
Topologically, the variety $\overline{J}_C$ is homeomorphic to a 2-dimensional sphere, and the perverse filtration is trivial by the Hard Lefschetz symmetry:
\[
P_0H^0(\overline{J}_C, \BQ) = H^0(\overline{J}_C, \BQ), \quad 0 = P_1H^2(\overline{J}_C, \BQ) \subset P_2H^2(\overline{J}_C, \BQ) = H^2(\overline{J}_C, \BQ).
\]

Next, we consider the case $(p,q) = (3,4)$. The curve $C$ has arithmetic genus $3$, and the compactified Jacobian $\overline{J}_C$ is a 3-fold. By \cite[Section 1.5]{MY} and the references therein, the variety $\overline{J}_C$ admits an affine paving
\[
\overline{J}_C = \BC^3 \sqcup \BC^2 \sqcup \BC^2 \sqcup \BC^1 \sqcup \BC^0; 
\]
furthermore,  we have
\[
\sum_{i=0}^6  \dim H^i(\overline{J}_C, \BQ) q^i= 1+q^2+2q^4 +q^6
\]
and
\[
\sum_{i=0}^6  \dim \mathrm{Gr}^P_i H^*(\overline{J}_C, \BQ) q^i= 1+q^2+q^3+ q^4 +q^6.
\]

In fact, for any $C$ associated with the singularity (\ref{p_q}), the cohomology ring and the perverse filtration for $\overline{J}_C$ are calculated explicitly in terms of tautological classes and their relations in \cite{OY, OY2} using representation theory techniques.

For $(p,q)=(3,4)$, the formula of \cite[Theorem 1.1.3]{OY2} gives the structure of the cohomology ring of $\overline{J}_C$:
\[
H^*(\overline{J}_C, \BQ) = \BQ[e_2, e_3]/(e_2e_3, e^2_3, e_2^4).
\]
Here $e_2 \in H^2(\overline{J}_C, \BQ), e_3 \in H^4(\overline{J}_C, \BQ)$ are the tautological generators given in \cite[Theorem 1.1.2]{OY2}. By \cite[Proof of Theorem 1.1.5]{OY2}, the perverse filtration can also be written down explicitly
\begin{align*}
    P_0=P_1 = \BQ 1 \subset P_2 = \BQ 1 \oplus \BQ e_2 \subset P_3 & = \BQ 1 \oplus \BQ e_2 \oplus \BQ e_3 \\ \subset P_4 = P_5 &  = \BQ1 \oplus \BQ e_2 \oplus \BQ e_3 \oplus \BQ e_2^2 \subset P_6 = H^*(\overline{J}_C, \BQ).
\end{align*}

\subsection{Independence of the family}

We outline a proof that the perverse filtration defined in Section \ref{sec2.1} is intrinsic, \emph{i.e.}, it only depends on the locally planar integral curve $C$. We follow \cite[Section 3.8]{MY}. The purpose is to emphasize that the fact that the perverse filtration is intrinsic is a consequence of the full support property given by Theorem \ref{thm1.1}.

We fix the curve $C$. Let $\CC \to B$ be an arbitrary family as in Section \ref{sec1} which realizes $C$ as a fiber $\CC_0$, from which we define the perverse filtration. We can also choose such a family to be a versal deformation $\CC_V \to V$ of the curve $C \to 0 \in V$. It suffices to show that the perverse filtrations on $H^*(\overline{J}_C, \BQ)$ defined via $B$ and $V$ respectively coincide, as this would allow us to compare the perverse filtration defined by any family with the one defined by a versal family.

Since the perverse filtration only depends on an analytic neighborhood of the point $0 \in B$, we can shrink $B$ so that there is a morphism $\iota: B \to V$ realizing $\CC \to B$ as the pullback of $\CC_V \to V$. Let 
\[
\pi: \overline{J}_\CC \to B, \quad \pi_V: \overline{J}_{\CC_V} \to V
\]
be the corresponding compactified Jacobian fibrations.

We pick isomorphisms
\[
R\pi_*\BQ_{\overline{J}_\CC} \simeq \bigoplus_i \CL_i[-i] \in D^b_c(B)
\]
and
\[
R\pi_{V*}\BQ_{\overline{J}_{\CC_V}} \simeq \bigoplus_i \CL'_i[-i] \in D^b_c(V)
\]
where $\CL_i, \CL'_i$ are semisimple perverse sheaves on $B,V$ respectively. By the base change, we have a canonical isomorphism
\[
 R\pi_*\BQ_{\overline{J}_\CC} = \iota^*R\pi_{V*}\BQ_{\overline{J}_{\CC_V}}  \in D^b_c(B)
\]
which yields isomorphisms
\begin{equation}\label{match1}
\bigoplus_{i} \CL_i[-i] \simeq     R\pi_*\BQ_{\overline{J}_\CC} = \iota^*R\pi_{V*}\BQ_{\overline{J}_{\CC_V}} \simeq \bigoplus_i \iota^*\CL'_i[-i] \in D^b_c(B).
\end{equation}
Here the first and the third isomorphisms are non-canonical. For our purpose, we need to show that the second canonical isomorphism yields an isomorphism of the perverse filtrations on $H^*(\overline{J}_C, \BQ)$; it suffices to prove
\[
\CL_i \simeq \iota^*\CL'_i \in \mathrm{Perv}(B).
\]

Since each term on the right-hand side of (\ref{match1}) is a semisimple perverse sheaf with full support by Theorem \ref{thm1.1}, each term on the left-hand side must also be a semisimple perverse sheaf with full support. Therefore, by comparing $\iota^*\CL'_i$ and $\CL_i$ over the locus $V^\circ$ of nonsingular curves, we conclude that they are isomorphic, which completes the proof. \qed

\begin{rmk}
    Let $f: X \to B$ be a proper morphism between nonsingular varieties with equi-dimensional fibers. Then the perverse $t$-structure yields a global perverse filtration
    \begin{equation}\label{fil1}
    P_0H^*(X, \BQ) \subset P_1H^*(X, \BQ) \subset P_2H^*(X, \BQ) \subset \cdots \subset H^*(X, \BQ)
    \end{equation}
    and a local perverse filtration for each fiber $X_b \subset X$,
    \begin{equation}\label{fil2}
     P_0H^*(X_b, \BQ) \subset P_1H^*(X_b, \BQ) \subset P_2H^*(X_b, \BQ) \subset \cdots \subset H^*(X_b, \BQ).
     \end{equation}
    Globally, (\ref{fil1}) is induced by taking the global cohomology of 
    \[
    \bigoplus_{i\leq k+ \dim B} H^*\left(B,{^\mathfrak{p}}\CH^i(Rf_* \BQ_X)[-i]\right) \hookrightarrow H^*(B, Rf_*\BQ) = H^*(X, \BQ)
    \]
    via a choice of isomorphism (\ref{DT_general}). Equivalently, it is given by taking the global cohomology of the perverse truncation functor
    \begin{equation}\label{truncation2}
    {^\mathfrak{p}}\tau_{\leq k+\dim B} Rf_* \BQ_{X} \to Rf_*\BQ_X.
    \end{equation}
    The local version is then induced by taking the stalk of (\ref{truncation2}) over $b \in B$. Both filtrations \emph{a priori} depend on $f: X \to B$.

    For the Hitchin fibration $h: M_{\mathrm{Higgs}} \to A$, the associated global perverse filtration plays a crucial role in non-abelian Hodge theory, as we will review in Section \ref{sec2.3}.
\end{rmk}

\subsection{A local P=W conjecture}\label{sec2.3}

A fascinating yet still mysterious feature of the perverse filtration for $\overline{J}_C$ is its conjectural connection to invariants associated with links. This serves as our motivation for exploring structures of the perverse filtration. This conjecture was proposed by Shende \cite{Sh0} aiming at unifying two stories as we briefly discuss below.

\subsubsection{Story 1: the P=W conjecture for the Hitchin fibration}\label{sec2.3.1}
Let $\Sigma$ be a compact Riemann surface of genus $g\geq 2$. Given two coprime integers $r,d$ (which will stand for the rank and the degree), there are two moduli spaces associated with such data: the moduli space $M_{\mathrm{Higgs}}$ of stable Higgs bundles of rank $r$ and degree $d$, and the corresponding character variety $M_{\mathrm{char}}$ parameterizing $d$-twisted rank $r$ complex representations of the fundamental group $\pi_1(\Sigma)$. These two moduli spaces are of different flavors. The variety $M_{\mathrm{Higgs}}$ underlies an algebraically complete integrable system, known as the Hitchin fibration
\[
h: M_{\mathrm{Higgs}} \to A \simeq \BC^{\frac{1}{2} \dim M_{\mathrm{Higgs}}};
\]
the character variety $M_{\mathrm{char}}$ is an affine variety only dependent on the topology of $\Sigma$. The non-abelian Hodge theory gives a diffeomorphism of these two moduli spaces. The $P=W$ conjecture of de Cataldo--Hausel--Migliorini \cite{dCHM}, proven recently \cite{MS_PW, HMMS, MSY}, says that the non-abelian Hodge theory exchanges the perverse filtration associated with the Hitchin fibration and the weight filtration associated with the mixed Hodge structure on $M_{\mathrm{char}}$:
\begin{equation}\label{P=W_conj}
``P=W": \quad P_kH^*(M_{\mathrm{Higgs}}, \BQ) = W_{2k} H^*(M_{\mathrm{char}}, \BQ).
\end{equation}
In particular, certain algebro-geometric invariants (\emph{e.g.} the invariants encoded in the perverse filtration) of the Hitchin fibration can be completely understood by the fundamental group $\pi_1(\Sigma)$ of the topological surface $\Sigma$.

\subsubsection{Story 2: link invariants associated with planar singularities}\label{sec2.3.2}

For a planar singular point $p \in C$, we can embed the neighborhood of the point $p$ into $\BC^2 \simeq \BR^4$. If we use a small 3-dimensional sphere $S^3_\epsilon \subset \BR^4$ to cut it, we obtain a link  
\[
L_{C,P}:= S^3_\epsilon \cap C \subset  S^3_{\epsilon}.
\]
A natural question is how the \emph{algebro-geometric} invariants of the curve interact with link invariants of $L_{C,P}$. This question was first proposed by Pandharipande using the Hilbert schemes (whose connection to the compactified Jacobian will be discussed in Section \ref{sec3}), and was answered conjecturally by Oblomkov--Shende \cite{OS}. Oblomkov--Shende's conjecture relates the invariants associated with the Hilbert schemes $C^{[n]}$ to the HOMFLY polynomial of the links $L_{C,P}$. This conjecture, as well as its colored refinement \cite{DHS}, was proven by Maulik \cite{Ma}. A further extension of the Oblomkov--Shende conjecture was proposed by Oblomkov--Rasmussen--Shende \cite{ORS}, which connects algebro-geometric invariants associated with a planar singularity with the \emph{Khovanov-Rozansky homology} \cite{KR} of the link. We assume that $C$ is a rational curve with a unique planar singular point $P \in C$.\footnote{Alternatively, we can work with the local compactified Jacobian associated with a planar singular point; we work with this global version only for convenience.} With this assumption, a version of the Oblomkov--Rasmussen--Shende conjecture is the following.

\begin{conj}[Oblomkov--Rasmussen--Shende \cite{ORS}]\label{Conj_ORS}
Let $L_{C,P}$ be the link associated with the singular point $P\in C$. Then we have
    \begin{equation}\label{ORS_eq}
\dim    \mathrm{Gr}^P_i H^d(\overline{J}_C, \BQ) = \mathrm{Coeff}_{q^{2(i-g_a(C))}t^d} \left[\mathrm{HHH}^{a=0}(L_{C,P}) \right].
    \end{equation}
\end{conj}
The Khovanov–Rozansky homology is triply graded, commonly described using the $a,q,t$-degrees; we refer to \cite{GKS} for more detail. The conjecture predicts that the perverse and the cohomological gradings for the compactified Jacobian are matched with the $q,t$-gradings of the lowest $a$-degree (\emph{i.e.} Hochschild degree) part, under a very explicit change of variables. We note that there is also a conjectural model for higher $a$-degrees using slightly more complicated moduli spaces \cite{ORS}.

\subsubsection{Shende's proposal}
As in Section \ref{sec2.3.2}, we assume that $C$ is rational with a unique planar singular point $P$. Shende proposed that Conjecture \ref{Conj_ORS} should be the numerical shadow of a ``local'' $P=W$ conjecture. The $P$-side concerns the algebraic geometry of the compactified Jacobian for a planar singularity, which is analogous to the Hitchin fibration. The $W$-side concerns a variety $M_{L_{C,P}}$ which should only be dependent on the link $L_{C,P}$, and is analogous to the character variety; see \cite{Trinh, Trinh2, BBAMY} for relevant discussions. Although Shende's proposal is still wide open, such a conjectural picture already imposes very strong constraints on the (well-defined) perverse filtration (\ref{P_fil}). We discuss them in Section \ref{decomp_?}.

\subsection{Predictions and known results}\label{decomp_?} 

First, the weight filtration for any algebraic variety is compatible with the cup-product; Shende's proposal predicts that the same holds for the perverse filtration (\ref{P_fil}). 

\medskip

\noindent { \bf Prediction A.} The perverse filtration (\ref{P_fil}) is compatible with the cup-product, \emph{i.e.}
\[
\cup: P_iH^d(\overline{J}_C, \BQ) \times P_{i'}H^{d'}(\overline{J}_C, \BQ) \to P_{i+i'} H^{d+d'}(\overline{J}_C, \BQ).
\]
\medskip

Second, recall that by a result of Shende \cite[Corollary 1]{Shende} the weight filtration $W_\bullet$ for the character variety $M_{\mathrm{char}}$ is of Hodge--Tate type; furthermore, what Shende proved is a stronger statement that the (decreasing) Hodge filtration $F^\bullet$ splits the (increasing) weight filtration over the rational coefficients:
\[
H^*(M_\mathrm{char}, \BQ) = \bigoplus_{k,d} \left(W_{2k} \cap F^k \right) H^d(M_{\mathrm{char}}, \BQ).
\]
This provides the weight filtration $W_\bullet H^*(M_{\mathrm{char}}, \BQ)$ a multiplicative splitting. If one believes that the cohomology of the local character variety $M_{L_{C,P}}$ associated with a link behaves like its global counter-part $H^*(M_{\mathrm{char}}, \BQ)$, then the following should hold.

\medskip

\noindent { \bf Prediction B.} The perverse filtration (\ref{P_fil}) admits a multiplicative splitting.

\medskip

We note that the perverse filtrations (\ref{fil1}) and (\ref{fil2}) \emph{a priori} are not multiplicative. This is due to the fact that the perverse $t$-structure is not compatible with the tensor product in the category $D^b_c(B)$. Understanding the multiplicativity was raised as a major challenge in proving the $P=W$ conjecture of Section \ref{sec2.3.1}; see \cite[Section 1]{dCHM} and \cite[Theorem 0.6]{dCMS}.

Prediction A was confirmed for the singularity (\ref{p_q}) by Oblomkov--Yun \cite{OY,OY2} using representation theory techniques. In fact, Oblomkov--Yun \cite{OY2} gave a complete description of the cohomology ring and the perverse filtration in terms of the tautological classes as we mentioned in Section \ref{Examples}. Therefore, in principle, it should be possible to check Prediction B. However, the calculation of the $\epsilon$-saturation in \cite[Theorem 1.1.3]{OY2} is too complicated to make such a verification practical. We leave this as a challenge to the interested reader.

%\begin{thm}[Oblomkov--Yun \cite{OY2}]
%    Both Predictions A and B hold for the rational curve with a singular point given by (\ref{p_q}).
%\end{thm}

In \cite{MSY}, using a different approach, Prediction A was confirmed completely.

\begin{thm}[Maulik--Shen--Yin \cite{MSY}]\label{thm2.6}
    Prediction A holds for any locally planar integral curve.
\end{thm}

The proof of Theorem \ref{thm2.6} relies on characterizing the perverse filtration (\ref{P_fil}) via Fourier transform which we will discuss in Section \ref{sec4}. We state Prediction B in the general case as a precise conjecture.

\begin{conj}\label{conj2.7}
    The perverse filtration (\ref{P_fil}) admits a multiplicative splitting. More precisely, there is a decomposition
    \[
    H^*(\overline{J}_C, \BQ)=\bigoplus_{i=0}^{2g} H^*_{(i)}(\overline{J}_C, \BQ)
    \]
    splitting the perverse filtration
    \[
    P_kH^*(\overline{J}_C, \BQ) = \bigoplus_{i\leq k} H_{(i)}^*(\overline{J}_C, \BQ),
    \]
    which is multiplicative with respect to the cup-product,
    \[
    \cup: H_{(i)}^*(\overline{J}_C, \BQ) \times H_{(i')}^*(\overline{J}_C, \BQ) \to H_{(i+i')}^*(\overline{J}_C, \BQ).     \]
\end{conj}

A subtle point regarding Conjecture \ref{conj2.7} is that we define the perverse filtration (\ref{P_fil}) via a fibration $\pi: \overline{J}_\CC \to B$ but a multiplicative splitting cannot exist globally for such a fibration. More precisely, it was shown in \cite{BMSY1} that for any $g\geq 3$, there exists a family $\CC \to B$ satisfying
\begin{enumerate}
    \item[(a)] each curve $\CC_b$ is of arithmetic genus $g$ with at worst one simple node, and 
    \item[(b)] the total space $\overline{J}_\CC$ is nonsingular,
\end{enumerate}
such that the perverse filtration $P_\bullet H^*(\overline{J}_\CC, \BQ)$ defined by (\ref{fil1}) does not admit a multiplicative splitting.\footnote{When $g=2$, such a family exists if one allows more simple nodes; see \cite[Proof of Theorem 0.4]{BMSY1}.} Consequently, for this family there is no sheaf-theoretic decomposition (\ref{DT_main}) compatible with the cup-product. An interesting feature of the examples in \cite[Theorem 0.4]{BMSY1} is that the cohomology of each fiber admits a multiplicative bigraded structure induced by a splitting of the perverse filtration, but these decompositions cannot be patched together.

New ideas may be needed to attack Conjecture \ref{conj2.7}.

By Theorem \ref{thm2.6}, taking the associated graded with respect to the perverse filtration yields a \emph{degeneration} of the graded $\BQ$-algebra
\[
\left(H^*(\overline{J}_C, \BQ):= \bigoplus_{d} H^d(\overline{J}_C, \BQ), ~~  \cup \right)
\]
to a \emph{bigraded} $\BQ$-algebra:
\[
\left(\BH_{C}: = \bigoplus_{i,d} \mathrm{Gr}^P_i H^d(\overline{J}_C, \BQ), ~~ \overline{\cup} \right)
\]
where the induced cup-product $\overline{\cup}$ on the associated graded is well-defined by Theorem \ref{thm2.6}. 

If Conjecture \ref{conj2.7} holds, then this degeneration is trivial; in particular, the cohomology $H^*(\overline{J}_C, \BQ)$ itself has a deep multiplicative \emph{hidden} bigrading.

\subsection{Further remarks}

\begin{enumerate}
    \item[(1)] Maulik--Yun constructed a Lefschetz (decreasing) filtration on $H^*(\overline{J}_C, \BQ)$ and conjectured in \cite[Conjecture 2.17]{MY} that this filtration splits the perverse filtration (\ref{P_fil}). This gives a natural candidate of a conjectural multiplicative splitting. In Section \ref{sec3} we will discuss another interesting splitting of the perverse filtration using the Hilbert schemes.
    \item[(2)] Theorem \ref{thm2.6} has been further extended to locally planar reduced curves. Given such a curve $C$, a model of the compactified Jacobian requires a choice of stability condition $\phi$; see \cite{MeRV}. For two choices of nondegenerate stability conditions $\phi, \phi'$ (that is, there is no strictly semistable locus), the corresponding fine compactified Jacobians $\overline{J}^\phi_C, \overline{J}^{\phi'}_C$ may be different topological spaces \cite[Theorem B (iii)]{MeRV}. However, it was proven in \cite[Theorem 0.8]{BMSY2} that their degenerations via the perverse filtrations are isomorphic bigraded $\BQ$-algebras:
    \[
    \left( \bigoplus_{i,d} \mathrm{Gr}^P_iH^d(\overline{J}^\phi_C, \BQ),~~\overline{\cup}^\phi \right)  \simeq  \left( \bigoplus_{i,d} \mathrm{Gr}^P_iH^d(\overline{J}^{\phi'}_C, \BQ),~~\overline{\cup}^{\phi'} \right).
    \]
    If the statement of Conjecture \ref{conj2.7} holds also for locally planar reduced curves, the result mentioned above implies (surprisingly) that the cohomology rings of $\overline{J}^\phi_C$ and $\overline{J}^{\phi'}_C$ are isomorphic on the nose. 
    \item[(3)] We refer to \cite{KT} for some recent progress on Conjecture \ref{Conj_ORS} via the geometry of Quot schemes.
    \item[(4)] From the perspective of the decomposition theorem for an abelian fibration, it was conjectured that the obstruction for having a \emph{multiplicative} decomposition is global; see \cite{Sh} and \cite[Section 0.6.3]{BMSY2}. This gives another explanation that the ``local'' perverse filtration (\ref{P_fil}) should admit a multiplicative splitting, which aligns with the expectation stated in Prediction B.     
    \end{enumerate}

\section{Hilbert schemes}\label{sec3}

\subsection{Compactified Jacobians vs Hilbert schemes}

The purpose of this section is to explain the connection between compactified Jacobians and Hilbert schemes of points. For a nonsingular curve, this is the classical Macdonald formula; for locally planar integral curves, this connection holds in a subtle way in which the perverse filtration comes into play naturally, as proven in the work of Maulik--Yun \cite{MY} and Migliorini--Shende \cite{MS}. By the result of Rennemo \cite{R}, Hilbert schemes also provide a geometric description of the perverse filtration as well as a natural splitting; we discuss this in Section \ref{sec3.3}.

Let $C$ be a locally planar integral curve of arithmetic genus $g$ with a smooth marked point $x \in C$. We consider the Hilbert scheme of points $C^{[k]}$. When $k \geq 2g-1$, the generalized Abel--Jacobi map 
\[
\mathrm{AJ}: C^{[k]} \to \overline{J}_C, \quad I_Z \mapsto I_Z \otimes \CO_C(kx)
\]
is a projective bundle. The geometry of $C^{[k]}$ is more interesting in the range $k < 2g-1$. Roughly, the cohomology group $H^d(C^{[k]}, \BQ)$ is responsible for the classes of $\mathrm{Gr}^P_kH^d(\overline{J}_C, \BQ)$.

\begin{thm}[\cite{MY, MS}]\label{thm3.1}
There is an isomorphism of graded $\BQ$-vector spaces
\begin{equation}\label{MCDN}
H^d(C^{[k]}, \BQ) \simeq \bigoplus_{i+j \leq k} \mathrm{Gr}_i^P H^{d-2j}(\overline{J}_C, \BQ).
\end{equation}
\end{thm}

Expanding (\ref{MCDN}), we see that
\[
H^d(C^{[k]}, \BQ) \simeq \mathrm{Gr}_k^PH^{d}(\overline{J}_C, \BQ) \oplus \mathrm{[lower~terms]}.
\]
Here, the lower terms contain those that have either a smaller perverse degree or a smaller cohomological degree.

\subsection{Theorem \ref{thm3.1} as a support theorem}

Theorem \ref{thm3.1} was proven by a support theorem analogous to Theorem \ref{thm1.1}. Here we sketch the proof of \cite{MY} which reduces Theorem \ref{thm3.1} to Theorem \ref{thm1.1}.

First, as in the proof of Theorem \ref{thm1.1}, we include $C$ as a fiber in a ``large enough'' family $p: \CC\to B$ with $C = \CC_0$ over $0 \in B$; for example, we may just take $\CC$ to be a versal deformation of $C$. We further assume that
\begin{enumerate}
    \item[(a)]  there is a section $s: B \rightarrow C$ which intersects with $\CC_0$ at a smooth point $x \in C = \CC_0$, and
    \item[(b)] there is a finite $B$-morphism
    \[
    \nu: \CC \to \mathbb{P}^1 \times B.
    \]
\end{enumerate} 
Indeed, since a versal family $\CC \to B$ containing $C =\CC_0$ as a fiber can be constructed using the Hilbert scheme of curves in a projective space $\BP^N$ \cite[Proposition 3.5]{MY}, if we shrink $B$, a general linear projection in $\BP^N$ yields a finite map $\CC_b \to \BP^1$ for each $b\in B$ which gives the map $\nu$ in (b). We further shrink $B$ to an \'etale neighborhood of $0\in B$ so that we can obtain a section as in (a).

We consider the relative Hilbert schemes of $k$ points:
\[
\pi_k: \CC^{[k]} \to B.
\]
The versal assumption ensures that the total space $\CC^{[k]}$ is nonsingular for any $k$; see \cite[Proposition 3.5]{MY}. As in Section \ref{sec1}, we use $B^\circ$ to denote the locus of nonsingular curves with $p^\circ: \CC^\circ \to B^\circ$ the restriction of $p: \CC \to B$.

We note that Theorem \ref{thm3.1} is a consequence of the following support theorem.

\begin{thm}\label{thm3.2}
    The decomposition of $R\pi_{k*} \BQ_{\CC^{[k]}}$ has full support for any $k\geq 1$.
\end{thm}

Indeed, the Macdonald formula shows that Theorem \ref{thm3.1} holds for a nonsingular curve $\CC_b \subset \CC$ over $b \in B^\circ$, where the perverse grading and the cohomological grading coincide in the sense that:
\[
\mathrm{Gr}^P_i H^*(J_{\CC_b}, \BQ) = H^i(J_{\CC_b}, \BQ).
\]
It is easy to turn the statement of Theorem \ref{thm3.1} for nonsingular curves into a sheaf-theoretic statement concerning local systems over $B^\circ$. If Theorem \ref{thm3.2} holds, then Theorem \ref{thm3.1} for $C = \CC_0$ can be obtained from taking the stalk over $0\in B$ for the intermediate extension of the relative version of Theorem \ref{thm3.1} over $p^\circ: \CC^\circ \to B^\circ$.

\begin{proof}[Proof sketch of Theorem \ref{thm3.2}]
The strategy of \cite{MY} is to use a descending induction on $k$. When $k \geq 2g-1$, the generalized Abel--Jacobi map $\mathrm{AJ}_b: \CC_b^{[k]} \to \overline{J}_{\CC_b}$ defined via the marked point $s(b) \in \CC_b$ is a projective bundle. Relatively over the base $B$, this gives a projective bundle
\[
\mathrm{AJ}_B: \CC^{[k]} \to \overline{J}_\CC
\]
defined via the section $s: B \to \CC $. Since the morphism $\pi_k:\CC^{[k]} \to B$ is the composition of $\mathrm{AJ}_B$ and $\pi: \overline{J}_\CC \to B$, we conclude for $k\geq 2g-1$ that $R\pi_{k*} \BQ_{\CC^{[k]}}$ has full support by Theorem \ref{thm1.1}.

To complete the proof, we need to show that any support of $\pi_{k}: \CC^{[k]} \to B$ also appears as a support of $\pi_{k+1}: \CC^{[k+1]} \to B$. Consider the incidence variety
\[
T_k \subset \left( \CC^{[k]} \times \BP^1 \right) \times_B \CC^{[k+1]}
\]
parameterizing $(Z, Z',a,b)$ where:
\begin{enumerate}
    \item[$\bullet$]  $b \in B$,
    \item[$\bullet$] $Z,Z' \subset \CC_b$ are of lengths $k, k+1$ respectively, \emph{i.e.} $[Z] \in \CC_b^{[k]}, [Z'] \in \CC_b^{[k+1]}$, such that $Z\subset Z'$,
    \item[$\bullet$] $a \in \BP^1$, 
    \item[$\bullet$] and the incidence condition is that the quotient of the ideal sheaves $I_{Z}/I_{Z'}$ is a skyscraper sheaf on $\CC_b$ over the point $a \in \BP^1$ via the finite map $\nu_b: \CC_b \to \BP^1$.
\end{enumerate} 

\medskip

Consider the natural projections
\begin{equation}\label{incidence}
    \begin{tikzcd}
T_k \arrow[r, ""] \arrow[d, ""]
& \CC^{[k+1]} \arrow[d, "\pi_{k+1}"]\\
\CC^{[k]} \times \BP^1 \arrow[r, "\pi'_k"] & B.
\end{tikzcd}
\end{equation}
Here $\pi'_k$ denotes the composition of the natural projection $\CC^{[k]} \times \BP^1 \to \CC^{[k]}$ and $\pi_k: \CC^{[k]} \to B$.
The incidence yields a relative correspondence over $B$,
\[
[T_k]: R\pi'_{k*} \BQ_{\CC^{[k]} \times \BP^1} \to R\pi_{k+1*} \BQ_{\CC^{[k+1]}}.
\]
A geometric study of the incidence (\ref{incidence}) shows that any simple summand of $R\pi'_{k*} \BQ_{\CC^{[k]} \times \BP^1}$ must appear as a simple summand of $ R\pi_{k+1*} \BQ_{\CC^{[k+1]}}$ via the correspondence $[T_k]$; see \cite[Page 10]{MY}. This completes the descending induction.
\end{proof}

\subsection{A natural splitting}\label{sec3.3}

Theorem \ref{thm3.1} describes the associated graded of the perverse filtration, as a $\BQ$-vector space, in terms of the cohomology of the Hilbert scheme. The work of Rennemo \cite{R} further describes the perverse filtration (\ref{P_fil}) explicitly via the Hilbert schemes $\CC^{[k]}$ and the generalized Abel--Jacobi map $\mathrm{AJ}: C^{[k]} \to \overline{J}_C$. 

For convenience, we work with \emph{homology}; the cohomological version can be recovered by the duality. Consider the bigraded vector space
\[
V(C): = \bigoplus_{k,i\geq 0} H_i(C^{[k]}, \BQ).
\]
The pushforward along the generalized Abel--Jacobi maps induces
\begin{equation}\label{AAAJJJ}
\mathrm{AJ}_*: V(C) \to H_*(\overline{J}_C, \BQ).
\end{equation}

Analogous to the operators constructed by Grojnowski and Nakajima for Hilbert schemes of points on nonsingular surfaces, Rennemo constructed the operators
\[
\mu_{+}([\mathrm{pt}]), ~ \mu_{-}([\mathrm{pt}]), ~   \mu_{+}([C]), \mu_{-}([C]) : V(C) \to V(C)
\]
using the incidence variety
\[
C^{[k,k+1]} \subset C^{[k]} \times C^{[k+1]},
\]
whose $(i,k)$-bidegrees on $V(C)$ are 
\[
(0,1), ~~ (-2,0),~~ (2,1),~~ (0,-1)
\]
respectively; he further showed that they satisfy the commutation relations
\[
[\mu_-([\mathrm{pt}]), \mu_+([C])] =  [\mu_-([C]), \mu_+([\mathrm{pt}])] = \mathrm{id}.
\]
Define   
\[
W:= \mathrm{Ker}\left(\mu_-([\mathrm{pt}]) \right) \cap \mathrm{Ker}\left( \mu_-([C]) \right) \subset V(C).
\]

\begin{thm}[Rennemo \cite{R}]\label{Thm3.3}
We have the following. 
 \begin{enumerate}
     \item[(a)] There is a natural isomorphism
     \[
     W\otimes \BQ[\mu_+([\mathrm{pt}]), \mu_+([C])] \xrightarrow{~\simeq~} V(C).
     \]
     \item[(b)] The generalized Abel--Jacobi map (\ref{AAAJJJ}) induces an isomorphism 
     \[
     \mathrm{AJ}_*: W \xrightarrow{~\simeq~} H_*(\overline{J}_C, \BQ).
     \]
 \end{enumerate}
\end{thm}

Since $V(C)$ is naturally bigraded, the primitive space $W \subset V(C)$ is also bigraded
\[
W = \bigoplus_{k,i} W_{k,i}, \quad W_{k,i}:= W\cap H_{i}(C^{[k]}, \BQ).
\]
Theorem \ref{Thm3.3} constructs a geometric bigraded decomposition 
\begin{equation}\label{decompp}
H_*(\overline{J}_C, \BQ)= \bigoplus_{k,i} D_kH_i(\overline{J}_C, \BQ),\quad D_kH_i(\overline{J}_C, \BQ):= \mathrm{AJ}_* \left(W_{k,i} \right).
\end{equation}
This decomposition turned out to describe the perverse filtration.

\begin{prop}[Rennemo \cite{R}]\label{prop3.4}
   Dualizing the decomposition (\ref{decompp}) yields a splitting of the perverse filtration (\ref{P_fil}).
\end{prop}

Rennemo asked in \cite[Question 1.4]{R} whether the (dual of the) $D$-decomposition of $H^*(\overline{J}_C, \BQ)$ is compatible with the cup-product; this construction provides a natural candidate for the desired splitting of Conjecture \ref{conj2.7}.

\subsection{Further remarks}
\begin{enumerate}
    \item[(1)] The correspondence between the Hilbert schemes and the perverse filtration for the compactified Jacobian is closely related to the correspondence between Pandharipande--Thomas invariants and the Gopakumar--Vafa invariants for Calabi--Yau 3-folds \cite{HST, MT}, where the perverse filtration appeared naturally in the definition of the Gopakumar--Vafa invariants. We refer to the paragraph after \cite[Corollary 2]{MS} for more discussions on this connection. Theorem \ref{thm3.1} was applied in  \cite{MT,DK} to prove interesting cases of this conjectural correspondence for curves with (3-fold!) singularities. 
    \item[(2)] Both Theorems \ref{thm3.1} and \ref{Thm3.3} were generalized to locally planar reduced curves. Theorem \ref{thm3.1} was generalized by Migliorini--Shende--Viviani \cite{MSV}. The Weyl algebra picture for the Hilbert scheme as in Theorem \ref{Thm3.3} was generalized by Kivinen \cite{K}. However, for the reducible cases, it is unclear how Kivinen's representation theoretic description is related to the perverse filtration as in Proposition \ref{prop3.4}.
    \item[(3)] For nonsingular curves lying in a del Pezzo surface, the relative Jacobian fibrations (with respect to certain degrees) admit compactifications known as the Le Potier moduli spaces \cite{LP}. Theorem \ref{thm3.1} has been applied in \cite{PSSZ} to prove cases of the $P=C$ conjecture, relating the perverse filtration to the Chern filtration associated with the Le Potier moduli space. See \cite{SZ} for further applications of these techniques for other surfaces.
\end{enumerate}

\section{Derived equivalences}\label{sec4}

\subsection{Perverse filtrations revisited}\label{sec4.1}

A challenge in studying the perverse filtration (\ref{P_fil}) is the abstract nature of the perverse truncation functor (\ref{truncation2}), which in general is difficult to compute. A geometric description of the perverse filtration (\ref{P_fil}) was introduced in Section \ref{sec3.3} using the Hilbert schemes. The purpose of this section is to introduce another description of (\ref{P_fil}) using the Arinkin--Fourier--Mukai transform of the derived category of \emph{coherent sheaves}. As an application, we explain in Section \ref{sec4.3} a proof of Theorem \ref{thm2.6} --- the perverse filtration associated with $H^*(\overline{J}_C, \BQ)$ is multiplicative.

The strategy is to generalize the ideas of Beauville \cite{B, B1} and Deninger--Murre \cite{DM} which characterized the Leray filtration of an abelian scheme using the Fourier--Mukai transform.

In the following, we first briefly recall the work of Beauville and Deninger--Murre in the special case of a Jacobian fibration. Let $p: \CC \to B$ be a smooth family of genus $g$ curves, and let $\pi: J_\CC \to B$ be the Jacobian fibration; it is an abelian scheme. We consider the decomposition theorem
\begin{equation}\label{DT}
R\pi_*\BQ_{J_C} \simeq \bigoplus_{i=0}^{2g} \wedge^i R^1p_* \BQ_\CC[-i]    
\end{equation}
where every term on the right-hand side is a (shifted) local system. The observation of \cite{B,B1, DM} is that the decomposition (\ref{DT}) is of a \emph{motivic} nature, where the algebraic cycles are produced by the Fourier--Mukai transform. To explain this, we recall the normalized Poincar\'e line bundle $\CP$ on the relative product $J_\CC \times_B J_\CC$ which induces the derived equivalences
\[
\mathrm{FM}_{\CP}: D^b \mathrm{Coh}({J}_\CC) \xrightarrow{\simeq}  D^b \mathrm{Coh}({J}_\CC), \quad \mathrm{FM}_{\CP^{-1}} : D^b \mathrm{Coh}({J}_\CC) \xrightarrow{\simeq}  D^b \mathrm{Coh}({J}_\CC)
\]
where we use $\CP^{-1}$ to denote the kernel of the inverse of the Fourier--Mukai transform $\mathrm{FM}_\CP$.

The normalized Poincar\'e line bundle $\CP$ induces a cycle class
\[
\mathfrak{F} = \sum_{i} \mathfrak{F}_i \in \mathrm{CH}^*(J_\CC \times_B J_\CC)_\BQ, \quad \mathfrak{F}_i:= \mathrm{ch}_i(\CP) \in \mathrm{CH}^i(J_\CC \times_B J_\CC)_\BQ,
\]
which gives a correspondence
\[
\mathfrak{F}: R\pi_* \BQ_{{J}_{\CC}} \to R\pi_* \BQ_{{J}_{\CC}}.
\]
The inverse $\CP^{-1}$ induces the inverse of the correspondence $\mathfrak{F}$:
\[
\mathfrak{F}^{-1} = \sum_{i} \mathfrak{F}^{-1}_i: R\pi_* \BQ_{{J}_{\CC}} \to R\pi_* \BQ_{{J}_{\CC}}
\]
which is the inverse of $\mathfrak{F}$. The operators $\mathfrak{F}, \mathfrak{F}^{-1}$ are regarded as the Fourier transforms.

We consider the correspondence given by 
\[
\mathfrak{q}_i:= \mathfrak{F}_i \circ \mathfrak{F}^{-1}_{2g-i} :   R\pi_* \BQ_{{J}_{\CC}} \to R\pi_* \BQ_{{J}_{\CC}}
\]
where $\circ$ is the composition of correspondences. Based on the ideas of Beauville \cite{B, B1}, Deninger--Murre \cite{DM} proves the following result.

\begin{thm}[Deninger--Murre \cite{DM}]
There is a canonical decomposition (\ref{DT}) such that the image of the correspondence 
\[
\mathfrak{q}_i:  R\pi_* \BQ_{{J}_{\CC}} \to R\pi_* \BQ_{{J}_{\CC}}
\]
recovers the $i$-th summand
\[
 \wedge^i R^1p_* \BQ_\CC[-i]  \subset R\pi_* \BQ_{J_\CC}.
\]
Consequently, the correspondence given by
\[
\mathfrak{l}_k:= \sum_{i\le k} \mathfrak{q}_i: R\pi_* \BQ_{J_\CC} \to R\pi_* \BQ_{J_\CC}
\]
recovers the classical truncation functor
\[
\tau_{\leq k} R\pi_* \BQ_{J_\CC} \to R\pi_* \BQ_{J_\CC}. 
\]
\end{thm}

 Now we consider the compactified Jacobian fibration $\pi: \overline{J}_\CC \to B$ associated with a family of locally planar integral curves $p: \CC \to B$ with $\overline{J}_\CC$ nonsingular. As we will introduce in Section \ref{sec4.2}, the normalized Poincar\'e line bundle for nonsingular curves admits an extension to a Cohen--Macaulay sheaf by Arinkin \cite{A2}, which induces Fourier transforms $\mathfrak{F}, \mathfrak{F}^{-1}$. Analogously to the smooth case, we consider the correspondence 
\[
\mathfrak{q}_i:= \mathfrak{F}_i \circ \mathfrak{F}^{-1}_{2g-i} :   R\pi_* \BQ_{\overline{J}_{\CC}} \to R\pi_* \BQ_{\overline{J}_{\CC}}.
\]

\begin{thm}[Maulik--Shen--Yin \cite{MSY}]\label{thm4.2}
    The operator
    \[
    \mathfrak{p}_k:= \sum_{i\leq k} \mathfrak{q}_i:   R\pi_* \BQ_{\overline{J}_{\CC}} \to R\pi_* \BQ_{\overline{J}_{\CC}}
    \]
    realizes the perverse truncation functor
    \[
    {^\mathfrak{p}}\tau_{\leq k+\dim B} R\pi_* \BQ_{\overline{J}_\CC} \to R\pi_* \BQ_{\overline{J}_\CC}.
    \]
\end{thm}

Theorem \ref{thm4.2} suggests that the complexity of the perverse truncation functor over the singular locus $B \setminus B^\circ$ is matched with the complexity of the Arinkin sheaf over $\overline{J}_\CC \setminus J_{\CC^\circ}$. We discuss more details about the Arinkin sheaf in Section \ref{sec4.2}; then in Section \ref{sec4.3} we sketch a proof for a sheaf-theoretic multiplicativity of the perverse truncation functor which implies Theorem \ref{thm2.6}.

\subsection{Arinkin--Fourier--Mukai transforms}\label{sec4.2}

As we saw in Section \ref{sec4.1}, a key geometric input in the characterization of the perverse truncation functor in Theorem \ref{thm4.2} is Arinkin's Fourier--Mukai kernel
\[
\CP \in \mathrm{Coh}\left( \overline{J}_\CC \times_B \overline{J}_\CC \right).
\]

Let
\[
\CP^\circ \in \mathrm{Pic}\left(J_{\CC^\circ} \times_{B^\circ} J_{\CC^\circ} \right)
\]
be the normalized Poincar\'e line bundle over the smooth locus. It is not hard to show that this bundle can be extended to a line bundle $\CP^{\heart}$ on the larger open subset
\begin{equation}\label{jjj}
 \left(J_\CC \times_B \overline{J}_\CC \right) \cup \left( \overline{J}_\CC \times_B {J}_\CC \right)
 \hookrightarrow \overline{J}_\CC \times_B \overline{J}_\CC
 \end{equation}
 In fact, over a pair of degree 0 line bundles
\[
(L_1, L_2) \in J_{\CC_b} \times J_{\CC_b}
\]
on a nonsingular curve $\CC_b$, the fiber of the normalized Poincar\'e line bundle is given by 
\[
\CP^\circ\Big{|}_{(L_1,L_2)} = \mathrm{det}\mathrm{R}\Gamma(\CC_b, L_1\otimes L_2) \otimes \mathrm{det}\mathrm{R}\Gamma(\CC_b, \CO_{\CC_b}) \otimes \mathrm{det}\mathrm{R}\Gamma(\CC_b, L_1)^{-1}\otimes \mathrm{det}\mathrm{R}\Gamma(\CC_b, L_2)^{-1}.
\]
This formula makes sense for singular curves $\CC_b$ as long as one of $L_1, L_2$ is still a line bundle; however, in the boundary where both factors are not line bundles (but only rank 1 torsion free sheaves), their (derived) tensor product may not be bounded. Therefore the first term in the right-hand side is not well-defined.

Nevertheless, if we use $j$ to denote the open embedding (\ref{jjj}), Arinkin showed in \cite{A2} that $\CP:=j_*\CP^{\heart}$ is a Cohen--Macaulay coherent sheaf on $\overline{J}_\CC \times_B \overline{J}_\CC$; he further proved that $\CP$ induces a derived auto-equivalence
\[
\mathrm{FM}_{\CP}: D^b\mathrm{Coh}(\overline{J}_\CC) \xrightarrow{~~\simeq~~} D^b\mathrm{Coh}(\overline{J}_\CC),
\]
whose inverse is induced by the kernel
\[
\CP^{-1} : = \CP^* \otimes \pi_2^* \omega_{\pi}[g], \quad \CP^*: = {\CH}om(\CP, \CO_{\overline{J}_\CC\times_B \overline{J}_\CC}).
\]
Here $\omega_\pi$ is the relative canonical bundle with respect to $\pi: \overline{J}_\CC \to B$, and $\pi_2: \overline{J}_\CC \times_B \overline{J}_\CC \to \overline{J}_\CC$ is the natural morphism.

Unlike the smooth case, the Fourier--Mukai kernels $\CP, \CP^{-1}$ are \emph{not} even perfect complexes. Therefore, defining the cycles of the Fourier transforms
\[
\mathfrak{F}, \mathfrak{F}^{-1} \in \mathrm{CH}_*\left(\overline{J}_\CC \times_B \overline{J}_\CC \right)_\BQ
\]
on the singular total space relies on  the Baum--Fulton--MacPherson map in the Riemann-Roch formula for singular varieties. Although it is very difficult to explicitly calculate $\mathfrak{F}, \mathfrak{F}^{-1}$ due to the singularities, one can prove the following quadratic ``Fourier vanishing'' relations 
\[
\mathfrak{F}^{-1}_i \circ \mathfrak{F}_j= 0, \quad i+j < 2g
\]
in the Chow group. Eventually Theorem \ref{thm4.2} was deduced from these relations. We note that a crucial geometric input in proving these relations is again the Severi inequality (\ref{severi}) --- a key player in the proof of the support theorem in Section \ref{sec1.3}!

\subsection{Multiplicativity via the Fourier transform}\label{sec4.3}

We conclude by explaining the proof of Theorem \ref{thm2.6}. In fact, it follows from a stronger sheaf-theoretic result.

\begin{thm}[Maulik--Shen--Yin \cite{MSY}]\label{thm4.3}
Let $p:\CC \to B$ and $\pi: \overline{J}_\CC \to B$ be as in Theorem \ref{thm4.2}. Then the cup-product
\[
\cup: R\pi_* \BQ_{\overline{J}_\CC} \otimes R\pi_* \BQ_{\overline{J}_\CC}  \to R\pi_* \BQ_{\overline{J}_\CC} 
\]
is compatible with 
the perverse truncation functor in the sense that
\[
\cup: {^\mathfrak{p}}\tau_{\leq k+\dim B} R\pi_* \BQ_{\overline{J}_\CC} \otimes{^\mathfrak{p}}\tau_{\leq k'+\dim B} R\pi_* \BQ_{\overline{J}_\CC} \to {^\mathfrak{p}}\tau_{\leq k+k'+\dim B} R\pi_* \BQ_{\overline{J}_\CC}.
\]
\end{thm}

It is clear that Theorem \ref{thm4.3} implies Theorem \ref{thm2.6}: if $C$ is a fiber $\CC_0$ over $0\in B$, we obtain Theorem \ref{thm2.6} by taking the stalk over $0\in B$.

The idea of the proof of Theorem \ref{thm4.3} is to use the interaction between the Fourier transform $\mathfrak{F}$ and the perverse truncation functor ${^\mathfrak{p}}\tau_{\leq \bullet} R\pi_* \BQ_{\overline{J}_\CC}$. Such a connection is a consequence of Theorem \ref{thm4.2}, as the perverse truncation functor admits an explicit formula in terms of the Fourier transforms. This connection leads to the observation that, in order to prove the multiplicativity, it suffices to control the size of the support of the convolution kernel 
\[
\CK \in D^b\mathrm{Coh}\left( \overline{J}_\CC \times_B \overline{J}_\CC \times_B \overline{J}_\CC \right);
\]
\emph{i.e.} $\CK$ is the Fourier--Mukai kernel for the convolution product
\[
\ast^\CK: D^b\mathrm{Coh}(\overline{J}_\CC) \otimes D^b\mathrm{Coh}(\overline{J}_\CC) \to   D^b\mathrm{Coh}(\overline{J}_\CC)
\]
which exchanges with the tensor product
\[
\otimes: D^b\mathrm{Coh}(\overline{J}_\CC) \otimes D^b\mathrm{Coh}(\overline{J}_\CC) \to   D^b\mathrm{Coh}(\overline{J}_\CC)
\]
under the Fourier--Mukai transforms $\mathrm{FM}_\CP$, $\mathrm{FM}_{\CP^{-1}}$.

Finally, Theorem \ref{thm4.3} follows from the fact that
\[
\mathrm{codim}_{\overline{J}_\CC\times_B \overline{J}_\CC \times_B \overline{J}_\CC}\left(\mathrm{supp}(\CK) \right) =g,
\]
whose proof relies (again!) on the Severi inequality (\ref{severi}).

\subsection{Further remarks}

\begin{enumerate}
    \item[(1)]  Theorem \ref{thm4.2} and its variants provide a tool to compute the perverse filtration for certain abelian fibrations. It was used in \cite{MSY} to give a proof of the $P=W$ conjecture, in \cite{MSY} to prove the Corti--Hanamura motivic decomposition conjecture for compactified Jacobian fibrations, in \cite{Ghosh} to locate the Chern classes in the perverse filtration, and in \cite{MSY,PSSZ} to prove the $P=C$ conjecture for del Pezzo surfaces for low degree cohomology.
    \item[(2)] The Arinkin sheaf, together with the induced Fourier--Mukai duality, was extended to locally planar reduced curves by Melo--Rapagnetta--Viviani \cite{MeRV2}. It is wide open whether it can be further extended to non-reduced locally planar curves. See \cite{Li} for some results in this direction. 
    \item[(3)] Although the Arinkin sheaf is not extended to the total space of the moduli of stable Higgs bundles, the Corti--Hanamura motivic decomposition conjecture was proven in \cite{MSY2} for the Hitchin fibration, including the locus of non-reduced spectral curves. The proof is based on a combination of the Arinkin sheaf and Springer theory.
\end{enumerate}

% -----------------------------------------------------------------------------
% Bibliography
% The bibliography is embedded in this .tex file, as requested.  No .bib file
% is needed.  Alphabetize entries by author surname and supply DOI/MR data when
% available.  The AMS will process the list with its Enhanced References system.
% Replace the sample entries rather than citing them in the final manuscript.
% -----------------------------------------------------------------------------

\end{document}